\documentclass[11pt]{amsart}

\usepackage{amsmath,amscd,amsfonts,amsthm,amssymb,eucal}
\usepackage{url}

  \newcommand{\shrinkmargins}[1]{
   \addtolength{\textheight}{#1\topmargin}
   \addtolength{\textheight}{#1\topmargin}
   \addtolength{\textwidth}{#1\oddsidemargin}
   \addtolength{\textwidth}{#1\evensidemargin}
   \addtolength{\topmargin}{-#1\topmargin}
   \addtolength{\oddsidemargin}{-#1\oddsidemargin}
   \addtolength{\evensidemargin}{-#1\evensidemargin}
   }

 \shrinkmargins{.5}

\renewcommand{\leq}{\leqslant}
\renewcommand{\geq}{\geqslant}

\DeclareMathOperator{\Hom}{Hom}

\DeclareMathOperator{\SL}{SL}

\DeclareMathOperator{\GL}{GL}

\DeclareMathOperator{\Spec}{Spec}

\DeclareMathOperator{\Sym}{Sym}

\DeclareMathOperator{\Gal}{Gal}

\DeclareMathOperator{\Conf}{Conf}

\DeclareMathOperator{\Aut}{Aut}
\DeclareMathOperator{\Sp}{Sp}

\DeclareMathOperator{\Br}{Br}

\DeclareMathOperator{\Epi}{Epi}
\DeclareMathOperator{\PSL}{PSL}
\DeclareMathOperator{\Out}{Out}

\newcommand{\field}[1]{\mathbb{#1}}

\newcommand{\Q}{\field{Q}}

\newcommand{\Z}{\field{Z}}
\newcommand{\Zhat}{\hat{\Z}}
\newcommand{\F}{\field{F}}

\newcommand{\R}{\field{R}}
\newcommand{\C}{\field{C}}
\renewcommand{\P}{\field{P}}
\newcommand{\A}{\field{A}}

\newcommand{\ra}{\to}

\newcommand{\HH}{\mathcal{H}}

\newcommand{\MM}{\mathcal{M}}

\newcommand{\set}[1]{\{#1\}}

\newcommand{\beq}{\begin{displaymath}}
\newcommand{\eeq}{\end{displaymath}}
\newcommand{\beqn}{\begin{equation}}
\newcommand{\eeqn}{\end{equation}}

\newcommand{\Fqbar}{\overline{\F}_q}

\newcounter{ccounter}

\theoremstyle{plain}
\newtheorem{thm}{Theorem}

\theoremstyle{definition}
\newtheorem{defn}[thm]{Definition}

\newtheorem{example}[thm]{Example}
\newtheorem{question}[thm]{Question}

\theoremstyle{remark}
\newtheorem{rem}{Remark}

\begin{document}
\title{Superstrong approximation for monodromy groups}
\author{Jordan S. Ellenberg}

\maketitle

This conference has been devoted to recent progress in the study of approximation properties for finitely generated subgroups of arithmetic lattices, especially such subgroups which are not themselves arithmetic.  Such groups are frequently encountered in geometry, where they appear as {\em monodromy groups} attached to the variation of cohomology in families of manifolds or algebraic varieties.

In the present paper we briefly survey the construction of monodromy groups, discuss our (limited) knowledge about whether such groups are arithmetic, and summarize the results of \cite{ehk}, which derive an application to arithmetic geometry from recent advances in superstrong approximation~\cite{helfgott},\cite{pyberszabo},\cite{breuillard-green-tao},\cite{golsefidy-varju}.  We conclude by indulging ourselves in some speculations about more general contexts, asking:  what are the interesting questions about ``nonabelian superstrong approximation" and ``superstrong approximation for Galois groups?"

\section*{Acknowledgments}

This document is an expanded version of a lecture presented by the author at a conference on ``Thin Groups and Superstrong Approximation" held at the Mathematical Sciences Research Institute in February 2012.  The work described here was partially supported by NSF Grant DMS-1101267 and a Romnes Faculty Fellowship.  The author is grateful to Ian Agol, Frank Calegari, and Peter Sarnak for useful conversations about the mathematics discussed here, and to MathOverflow for useful input concerning the image of the point-pushing subgroup in the symplectic group and the meaning of superstrong approximation for Galois groups.

\section{What is monodromy?}

The simplest example of monodromy is provided by the topological notion of a {\em covering space}.  Suppose given a path-connected base space $B$, endowed with a choice of basepoint $b$, and let $f:  X \ra B$ be a covering space of degree $n$; that is, each point $x \in X$ has a neighborhood $U$ such that $f^{-1}(U)$ is homeomorphic to $n$ disjoint copies of $U$.  In particular, $f^{-1}(b)$ consists of $n$ points.  A paradigmatic example is provided by $B = X =  \C^\times$, where the map $f$ is given by $z \mapsto z^n$, and $b=1$.

Now let $\gamma$ be a loop in $B$ beginning and ending at $b$, and let $\tilde{b}$ be a point in $f^{-1}(b)$.  It follows from the definition of covering space that there is a unique path $\tilde{\gamma}$ which starts at $\tilde{b}$ and projects to $\gamma$.  Intuitively, one imagines the point in $X$ ``following its shadow" in $B$, whose motion is specified by $\gamma$.  In the paradigmatic example, we can take $\gamma:[0,1] \ra B$ to be a counterclockwise traversal of the unit circle
\beq
\gamma(t) = e^{2 \pi i t}
\eeq
and $\tilde{b} = 1$.  Then $\tilde{b}$ starts at $1$, and at time $t$ is required to satisfy
\beq
\tilde{\gamma}(t)^n = e^{2 \pi i t}
\eeq
It follows that $\tilde{\gamma}(t) = e^{2 \pi i t / n}$.  The path $\tilde{\gamma}$ starts at $1$, but it ends at $e^{2 \pi i / n}$.  In fact, starting at any point $\tilde{b}$ and ``following your shadow" around $\gamma$ puts you at $\tilde{b} e^{2 \pi / n}$.  The permutation of the $n$ points in $f^{-1}(b)$ induced by this operation is called the {\em monodromy} of the cover around $\gamma$; in this case, it is a cyclic permutation of the $n$ points.  It is easy to check that the map sending a path $\gamma$ based at $b$ to its monodromy is a homomorphism
\beq
\mu:  \pi_1(B,b) \ra \Sym(f^{-1}(b))
\eeq
which we call the monodromy map of the cover.  The image of $\mu$ is called the {\em monodromy group}.  

Any identification of $f^{-1}(b)$ with the set $[1,\ldots,n]$ allows us to express the monodromy of the cover as a homomorphism from the fundamental group $\pi_1(B,b)$ to $S_n$; of course, changing the labeling on $f^{-1}(b)$ changes the monodromy by conjugation.  Note that the monodromy map need not surject onto the symmetric group of $f^{-1}(b)$; for instance, in the paradigmatic example, the monodromy group is cyclic of order $n$.  For that matter, the monodromy group can be trivial -- this happens precisely when $X$ is isomorphic to the trivial cover $B \times [1\ldots,n]$.  The monodromy group is transitive if and only if $X$ is connected; in this case, there is a well-defined conjugacy class of point stabilizers in the monodromy group, and the preimage in $\pi_1(B,b)$ of such a stabilizer is an index-$n$ subgroup of $\pi_1(B,b)$.  In this setting, we are just rephrasing the usual Galois-theoretic identification between connected degree-$n$ covers of $B$ and index-$n$ subgroups of $\pi_1(B,b)$.

\medskip

We can generalize this story:  suppose now that $f: X \ra B$ is a {\em fibration} of manifolds.  Denote the fiber over $b \in B$ by $X_b$.   One can take hold of a homology class in $H_i(X_b,\Z)$ and ``move it along $\gamma$" just as we did in above with a single point of $X$, and this defines a monodromy map
\beq
\mu:  \pi_1(B,b) \ra \GL(H_i(X_b,\Z))
\eeq

Again, the {\em monodromy group} of the fibration (or ``the monodromy group on $H_i$", if there is some ambiguity) is the image of $\pi_1(B,b)$ in the arithmetic group $\GL(H_i(X_b,\Z))$.  Of course, we can define a monodromy group using the action on cohomology in just the same way.

Rather than define this process rigorously, we will sketch an illustrative case.  Suppose that $B$ is the moduli space $\MM$ of lattices in $\R^2$ up to homothety, and let $X$ be the universal torus over $\MM$; in other words, for each point $p$ of $\MM$ parametrizing a lattice $\Lambda_p \subset \R^2$, the fiber $X_p$ is the torus $\R^2 / \Lambda_p$.

Now let the basepoint $b$ be the point parametrizing the square lattice $\Lambda_b$ generated by $(1,0)$ and $(0,1)$, and define a loop $\gamma \subset \MM$ by letting $\gamma(t)$ be the lattice generated by $(1,0)$ and $(t,1)$.  Note that $\gamma(1) = \gamma(0) = b$.

The path from $(0,0)$ to $(0,1)$ represents a nontrivial homology class in $H_1(X_b,\Z)$, which we call $e_1$.  Write $e_2$ for the homology class of the path from $(0,0)$ to $(1,0)$; then $e_1$ and $e_2$ freely generate $H_1(X_b,\Z)$.

What happens as we move along $\gamma$, carrying $e_1$ along with us as we go?  At time $t$, our path moves along a straight line segment from $(0,0)$ to $(t,1)$.  And when the loop is finished, our path goes from $(0,0)$ to $(1,1)$.  The torus has returned to its starting point, but the homology class on $X_b$ represented by our path is now $e_1 + e_2$.

What if we'd started with the path from $(0,0)$ to $(1,0)$?  This path stays in place through our whole trip around $\gamma$.

In other words, the loop $\gamma$ induces a {\em unipotent} monodromy transformation on $H_1(X_b,\Z)$, sending $e_1$ to $e_1 + e_2$ and fixing $e_2$.

\medskip

Below we present a few examples of monodromy maps attached to families of manifolds.

\begin{example}
If $X \ra B$ is a finite covering space of degree $n$, then $X_b$ is a set of $n$ points, and the only nontrivial homology group of $X_b$ is $H_0(X_b) \cong \Z^n$.  In this case, after choosing an identification of $X_b$ with $1,\ldots,n$,  the monodromy map factors as
\beq
\mu: \pi_1(B,b) \ra S_n \ra \GL_n(\Z)
\eeq
where the second map is the standard permutation representation.
\end{example}

\begin{example}
Let $B$ be the moduli space of genus $g$ Riemann surfaces, and let $X$ be the universal curve over $B$; that is, for each point $p$ of $B$, the fiber $X_p$ is the Riemann surface parametrized by $p$.  In this case, $\pi_1(B,b)$ is the {\em mapping class group} $\Gamma_g$ of genus $g$.  (See \cite{farbmargalit} for a thorough introduction to the mapping class group and its basic properties.)  We can identify $H_1(X_b,\Z)$ with $\Z^{2g}$; having done so, the monodromy map takes the form
\beq
\mu: \Gamma_g \ra \GL_{2g}(\Z)
\eeq
This map is definitely not surjective once $d > 2$, because $H_1(X_b,\Z)$ carries a symplectic intersection form which is preserved when moving around in moduli space.  So the image lies in the subgroup of matrices preserving this form, and we can write:
\beq
\mu:  \Gamma_g \ra \Sp_{2g}(\Z).
\eeq
In fact, the image of $\mu$ is known to be all of $\Sp_{2g}(\Z)$.  (See e.g. \S 6.3.2 of \cite{farb:farbmargalit})
 \end{example}

\begin{example}
Let $B$ be the configuration space parametrizing unordered $n$-tuples of distinct elements of $\C$, so that $\pi_1(B,b)$ is the Artin braid group $B_n$ on $n$ strands.  Let $X$ be the family whose fiber over a configuration $\set{p_1, \ldots, p_n}$  is the smooth proper algebraic curve with affine model
\beq
y^d = (x-p_1) (x-p_2) \ldots (x-p_n).
\eeq
Then the monodromy map
\beq
\mu:  B_n \ra \GL(H_1(X_b))
\eeq
is a representation of the braid group, which is studied extensively by McMullen in \cite{mcmu:hodge}.  Each fiber $X_b$ is an algebraic curve with a $\Z/d\Z$-action:  namely, $1$ sends $y$ to $e^{2 \pi i / d} y$ and leaves $x$ unchanged.  Thanks to the intersection form, the image of $\mu$ is again contained in the group of symplectomorphisms $H_1(X_b)$.  By contrast with the previous example, the monodromy group is not the full symplectic group, because it commutes with the canonical action of $\Z/d\Z$ on the fiber.  So the monodromy group is contained in the centralizer of $\Z/d\Z$ in the symplectic group; we denote this group by $\Sp_T(H_1(X_b))$.   
\label{ex:mcmullen} 
\end{example}

\begin{rem}  Another classical route to monodromy groups is via solutions of differential equations.  We won't pursue examples of this kind here, but see the papers in this volume by Sarnak~\cite{sarnak:msri} and Fuchs~\cite{fuchs:msri} for a fuller account.

\end{rem}

\section{Big monodromy and superstrong approximation}
\label{s:superstrong}

Given a fibration $X \ra B$, we might aim to describe, as best we can, the corresponding monodromy group $\Gamma = \mu(\pi_1(B,b))$  in $\GL(H_i(X_b,\Z))$.   A crucial invariant is the Zariski closure $G$ of $\Gamma$, which is an algebraic group over $\Q$ contained in a general linear group $\GL_n$, where $n = \dim_\Q H_i(X_b,\Q)$.  We refer to $G$ as the {\em algebraic monodromy group}.  By $G(\Z)$ we shall mean the intersection of $G(\C)$ with $\GL(H_i(X_b,\Z))$.\footnote{In this article, we will always be in a context where $H_i(X_b,\Z)$ is torsion-free, and where any subtleties about the ``right" model for $G$ over $\Spec \Z$ do not intervene.}

As a first step, we would like to be able to say whether it is {\em big}.  Of course, there are several senses in which the word ``big" could be meant.   For instance, we can ask whether $G$ is all of $\GL_n$.  As we've seen, the answer is often negative, sometimes for readily perceived reasons -- for instance, in the case where $X \ra B$ is a family of genus $g$ curves, the intersection form on $H_1(X_b,\Z)$ constrains $G$ to be contained in $\Sp_{2g}$.  So one might also ask, more vaguely, ``is $G$ as big as it could be, given all the reasons I can think of for it to be small?"  For a typical family of curves, this amounts to asking whether $G$ is the full symplectic group.  In the context of Example~\ref{ex:mcmullen}, it would mean to ask whether $G$ is the whole of $\Sp_T$.

It is not straightforward to compute the algebraic monodromy group $G$ in general, but much progress has been made towards ensuring that $G$ is big under various checkable conditions.  (See, for instance, Katz~\cite{katz2004larsen} and Hall~\cite{hall2008big} for representative recent results.) 

\subsection{Thin monodromy}

Even if $G$ is known, there are still important questions about the size of $\Gamma$ itself.  First of all, we can compare $\Gamma$ with the full lattice $G(\Z)$.  The critical question is whether $\Gamma$ has finite index in $G(\Z)$.

\begin{defn} A subgroup $\Gamma$ of $G(\Z)$ is {\em thin} if it is Zariski-dense and has infinite index.
\end{defn}

There is, of course, a large literature on the geometry, group theory, and dynamics of arithmetic lattices -- the main theme of this conference is to understand the extent to which some of these results extend to thin subgroups.

It can be difficult, in practice, to know whether the monodromy group of a family is thin, and both thin and non-thin cases do arise.  For instance, the family of hyperelliptic curves of genus $g$ has monodromy group which is finite-index in $\Sp_{2g}(\Z)$, by a theorem of A'Campo~\cite{acampo}.  On the other hand, the family of curves described in Example~\ref{ex:mcmullen} can have either thin or non-thin monodromy depending on the parameters;  for instance, McMullen~\cite[\S 11]{mcmu:hodge} shows (following Deligne and Mostow) that the monodromy is thin when $(n,d) = (4,18)$.  On the other hand, Venkataramana~\cite{venkataramana} shows that the monodromy is an arithmetic lattice whenever $n \geq 2d$.  Brav and Thomas~\cite{brav} have recently shown that a family of quintic three-folds constructed by Dwork has as its $H^3$-monodromy a thin subgroup of $\Sp_4(\Z)$.  Fuchs's article in this volume~\cite{fuchs:msri} describes several senses in which thinness seems to be the ``generic" situation, including many examples of differential equations whose monodromy groups appear on experimental grounds to be thin.

\begin{example}  Consider the family of complex algebraic curves given by 
\beq
X_t:  y^2 = f(x)(x-t)
\eeq
where $f(x)$ is a fixed squarefree polynomial of degree $2g$, with $g > 1$, and $t$ varies over the base $B$, which is the complement of the roots of $f$ in $\A^1$.  The monodromy group $\Gamma$ of this family is Zariski dense in $\Sp_{2g}$; see unpublished work of Yu~\cite{yu} and the more general theorem of Hall~\cite[Th 4.1]{hall2008big}.  In fact, Yu proves that the monodromy group of this family is a specified finite-index subgroup of $\Sp_{2g}(\Z)$, contained in the principal congruence subgroup of level $2$ and containing the principal congruence subgroup of level $4$.  We take this opportunity to correct an error in \S 4.1 of \cite{ehk}, where we make the incorrect assertion that the monodromy group is the whole of $\Gamma(2)$.  (This error does not affect the argument in which it appears.)

Yu computes the monodromy group by direct computation with matrices.  Ian Agol explained to us a very handsome pure-thought argument via the Margulis normal subgroup theorem; we include this argument here in the hope that the idea might be useful in other similar cases.

Let $\Conf_{2g+1}$ be the configuration space parametrizing squarefree complex polynomials of degree $2g+1$.  The fundamental group of $\Conf_{2g+1}$ is the Artin braid group $\Br_{2g+1}$, and there is a natural family of genus-$g$ curves over $\Conf_{2g+1}$ whose fiber over a polynomial $g(x)$ is the hyperelliptic curve $y^2=g(x)$.  By the theorem of A'Campo already mentioned, the image of the corresponding monodromy map
\beq
\Br_{2g+1} \ra \Sp_{2g}(\Z)
\eeq
is a finite-index subgroup of $\Sp_{2g}(\Z)$ containing the full level-2 congruence subgroup $\Gamma(2)$.

The permutation action on strands induces a map $\Br_{2g+1} \ra S_{2g+1}$ whose kernel $P_{2g+1}$ is called the {\em pure braid group}.  The image of $P_{2g+1}$ under the monodromy image above is precisely $\Gamma(2)$.

Now the base space $B$ considered in the present example maps to $\Conf_{2g+1}$ by sending $t$ to $f(x)(x-t)$, and the monodromy map $\pi_1(B,b) \ra \Sp_{2g}(\Z)$ thus factors as a composition
\beq
\mu:  \pi_1(B,b) \ra \pi_1(\Conf_{2g+1},b) \cong \Br_{2g+1} \ra \Sp_{2g}(\Z).
\eeq

Topologically speaking, the image of $\pi_1(B,b)$ in $\Br_{2g+1}$ is the {\em point-pushing group} $H$ consisting of all braids in which the first $2g$ strands are fixed in place, and the last is free to wind around the others.  The point-pushing group is the kernel in the {\em Birman exact sequence}
\beq
 H \ra P_{2g+1} \ra P_{2g} \ra 1
\eeq
and in particular is a {\em normal} subgroup of the pure braid group.  Thus, $\mu(H)$ is a normal subgroup of the lattice $\Gamma(2)$.  It then follows form the Margulis normal subgroup theorem that $\mu(H)$ is either finite or finite-index; but it is easy to check that a loop carrying $t$ once around a root of $f$ has unipotent (and thus infinite-order) monodromy, so $\mu(H)$ must be finite index, as claimed.

\label{ex:pointpush}
\end{example}

\bigskip

To sum up:  the monodromy group $\Gamma$ can be ``big" in the sense that its Zariski closure $G$ is as large as possible, but it might well be the case that $\Gamma$ is very ``small" by virtue of being thin, and in a given explicit geometric situation we may not know which is the case.  The question is then whether thin groups are ``big enough" for applications.

\subsection{Strong approximation}

For instance, one may ask how well $\Gamma$ approximates $G(\Z)$ locally -- that is, with respect to the $p$-adic Lie groups $G(\Z_p)$.  A natural question in this vein is whether the inclusion $\Gamma \ra G(\Z_p)$ is dense for almost all $p$; that is, does $\Gamma$ satisfy {\em strong approximation?}  Of course, this requires that $G(\Z)$ itself satisfies strong approximation; this is the case whenever $G$ is semisimple and simply connected and $G(\R)$ is noncompact.  It is a remarkable fact, due in its strongest form to Weisfeiler, that (under mild conditions on $G$) {\em any} finitely-generated Zariski-dense subgroup $\Gamma$ of $G(\Z)$ satisfies strong approximation.  For details and history concerning these results, see Rapinchuk's survey in this volume~\cite{rapinchuk:msri}.  

A central theme of the lectures in this workshop has been that, despite the name ``strong approximation," density of $\Gamma$ in almost all $G(\Z_p)$ is not the strongest local approximation condition we can place on $\Gamma$.  Indeed, strong approximation is {\em not} strong enough to make sieving arguments work.  For this we need a stronger condition, called {\em superstrong approximation}, which we explain below.  See the articles by Golsefidy~\cite{golsefidy:msri} and Kowalski~\cite{kowalski:msri} in this volume for a survey of the ``affine sieve" and other arguments which use superstrong approximation to derive quantitative statements on thin subgroups of arithmetic groups.

\subsection{Superstrong approximation}

Let $\Gamma$ be a finitely generated subgroup of $G(\Z)$, and let $\underline{\gamma} = \gamma_1, \ldots, \gamma_r$ be a set of generators for $\Gamma$.  If $\pi: G(\Z) \ra Q$ is a finite quotient of $G(\Z)$, then we can form a Cayley graph $X(Q,\underline{\gamma})$, whose vertices are labeled by the elements of $Q$, and which has $q$ adjacent to $q'$ precisely when $q' = q \pi(\gamma_i)^{\pm 1}$ for some generator $\gamma_i$.

In what follows, we will let $Q$ range over the congruence quotients $G(\Z/p\Z)$, and we will refer to the Cayley graph $X(G(\Z/p\Z),\underline{\gamma})$ simply as $X_p$.  

Note first that $\Gamma$ surjects on $G(Z/p\Z)$ if and only if $X_p$ is connected.  The notion of connectedness, in turn, can be addressed spectrally:  the (suitably normalized) Laplacian operator on $X_p$ always has an eigenvalue of $0$, attached to the space of constant eigenfunctions, and all other eigenvalues $\geq 0$.  In fact, the number of components of $X_p$ is precisely the multiplicity of the eigenvalue $0$.  In other words, the second-largest eigenvalue $\lambda_1$ is $0$ precisely when $X_p$ is not connected.  And strong approximation, from a spectral point of view, says that $\lambda_1(X_p) = 0$ for only finitely many $p$.

Replacing this condition with a quantitative bound yields the notion of {\em superstrong approximation}.

\begin{defn}  Let $\Gamma$ be a finitely generated subgroup of $G(\Z)$ and define a sequence of Cayley graphs $X_p$ as above.  We say $\Gamma$ has {\em superstrong approximation} if there exists a constant $\epsilon>0$ such that $\lambda_1 > \epsilon$ for all sufficiently large primes $p$.
\end{defn}

In other words, superstrong approximation is the condition that the $X_p$ form a family of {\em expander graphs}.  By the usual theory of expander graphs, this condition does not depend on the choice of the generating set $\underline{\gamma}$;  in other words, it is a property of $\Gamma \subset G(\Z)$, as our phrasing suggests.  We note that superstrong approximation follows immediately if $\Gamma$ has Kazhdan's property $T$, as it must if it is a lattice in a higher-rank arithmetic group.  But when $\Gamma$ is thin, superstrong approximation doesn't follow from any abstract group-theoretic property of $\Gamma$.  Indeed, in practice $\Gamma$ is often free, so has many actions on finite sets whose corresponding Cayley graphs have small $\lambda_1$; superstrong approximation demands that the $G(\Z/p\Z)$ are not among those finite sets.

The remarkable progress of the last few years towards uniform bounds for expansion in Cayley graphs of linear groups, as described elsewhere in this volume~\cite{pyber:msri},  has brought us to the point where we now have a superstrong analogue of Weisfeiler's theorem; under mild conditions on the ambient algebraic group $G$, thin groups have superstrong approximation merely by virtue of their Zariski density in $G$.

\section{Application:  gonality growth in families of covers}

\label{s:ehk}

We now describe an arithmetic application of superstrong approximation for monodromy groups.  All the material here is explained in substantially more detail in \cite{ehk}.

Let $U$ be an smooth (but not necessarily proper) algebraic curve over a number field $k$, and let $A \ra U$ be an abelian scheme over $U$ (i.e. a family of abelian varieties parametrized by $U$.)  If $t$ is a point of $U(\bar{k})$, we denote by $A_t$ the fiber of the family over $t$; it is an abelian variety defined over the finite extension $k(t)/k$.  If $n$ is an integer, $A_t[n]$ denotes the group of $n$-torsion points of $A_t$.

For instance, we could fix a squarefree polynomial $f$, let $U$ be the complement in $\A^1$ of the roots of $f$, and take $A$ to be the family whose fiber $A_t$ is the Jacobian of the hyperelliptic curve $y^2 = f(x)(x-t)$.  (Though be warned that to prove Theorem~\ref{th:sample} in this case doesn't use new results on thin groups, since in this case the monodromy of the family is finite-index in $\Sp_{2g}(\Z)$, as explained in Example~\ref{ex:pointpush}.)
 
The following theorem follows from Theorem 7 of \cite{ehk}.

\begin{thm}[Theorem 7, \cite{ehk}] For every $d \geq 1$ there exists a constant $\ell(d)$ such that, for all primes $\ell > \ell(d)$, there are only finitely many $t \in U(\bar{k})$ such that $[k(t):k] \leq d$ and $A_t[\ell](k(t))$ contains a nontrivial point.
\label{th:sample}
\end{thm}

In words, ``there are only finitely many $\ell$-torsion points on fibers of $A$ over points of bounded degree."

We now sketch the proof.  The first step is to construct finite covers of $U$ as follows.  For each integer $n$, let $f_n: U_n \ra U$ be the map whose fiber is the finite set of nontrivial $n$-torsion points $A_t[n]$; this is a finite unramified covering of degree $n^{2\dim A} - 1$.   One might think of these covers as analogues of the classical modular curves; for if $U$ is the moduli space of elliptic curves $Y(1)$\footnote{Warning:  $Y(1)$ is not really an algebraic curve but a stack with nontrivial generic inertia.  So it would be more precise to say that when $U = Y_1(p)$ for $p > 3$ and not dividing $n$, we have $U_n = Y_1(np)$.}, the cover $U_n$ is precisely the moduli space $Y_1(n)$ of elliptic curves with an $n$-torsion point.  In fact, the argument below is very much in the spirit of those of Abramovich~\cite{abramovich} and Zograf~\cite{zograf} which gave a lower bound for the gonality of the classical modular curves; where we use new developments in superstrong approximation, they use Selberg's theorem that (in modern language) $\SL_2(\Z)$ has property $\tau$ with respect to the family of congruence subgroups.

Theorem~\ref{th:sample} is really a statement about algebraic points on the covers $U_\ell$:  we are saying that, for all $\ell > \ell(d)$, the curve $U_\ell$ has only finitely many points whose field of definition has degree at most $d$ over $k$.

This statement might at first seem rather strong, since there are infinitely many degree-$d$ field extensions of $k$.  Even if a curve $C$ has genus at least $2$, so that $|C(k')| < \infty$ for each $k'/k$ of degree $d$ by Faltings' theorem, it is far from clear that the union of $C(k')$ over {\em all} such $k'$ should be finite.  Indeed, this is often not the case:  for instance, the curve
\beq
y^2 = x^5 + x + 1
\eeq
has infinitely many points over quadratic extensions of $\Q$, since one may produce such a point by specifying $x \in \Q$ arbitrarily.  Of course, the story is the same for any hyperelliptic curve, which may have arbitrarily large genus.  So the genus is not sufficient to distinguish curves with many degree-$d$ algebraic points from curves with few such points.  For that purpose, we need to study a different geometric invariant:  the {\em gonality}.

\begin{defn}  The {\em gonality} of an algebraic curve $C/k$ is the minimum integer $\gamma$ such that there exists a morphism $C \ra \P^1$ of degree $\gamma$.
\end{defn}

(In this section, we take gonality to be a {\em geometric} invariant of the curve; that is, in the definition we allow the morphism $C \ra \P^1$ to be defined over an extension of $k$.)

For example, a curve has gonality $1$ precisely when it has genus $0$, and gonality $2$ precisely when it is hyperelliptic but not rational.

The gonality of a curve has a very strong effect on its points over fields of bounded degree over a fixed field, thanks to work of Faltings~\cite{faltings}, Abramovich-Voloch~\cite{abramovichvoloch}, and Frey~\cite{frey}.  Namely, if $C/k$ has gonality greater than $2d$,  then $C(\bar{k})$ has only finitely many points $t$ such that $[k(t):k] \leq d$.

To prove Theorem~\ref{th:sample}, it thus suffices to show that the gonality $\gamma(U_\ell)$ increases with $\ell$.  This is accomplished using two main tools.  First, we pass to the complex numbers and invoke an analytic theorem of Li and Yau~\cite{liyau} which provides lower bounds for the gonality of a complex algebraic curve  in terms of the spectrum of the Laplacian.  In particular, it is enough to show that the family of Riemann surfaces $U_\ell$ has a Laplacian spectral gap.

Why should this be so?  This requires a combinatorial detour.   The argument runs roughly as follows.  Choose a basepoint $u \in U(\C)$ and a basis $\gamma_1, \ldots, \gamma_r$ for $\pi_1(U(\C),u)$.  (For simplicity we assume here that $U$ is a punctured sphere, so that $\pi_1(U)$ is free.)  We can represent each $\gamma_i$ by a loop in $U(\C)$ based at $u$, and we can arrange our loops to be mutually disjoint.  We denote the union of these loops by $X$; it is a one-vertex graph embedded in the Riemann surface $U(\C)$.

Now the preimage of $X$ under the cover $U_\ell \ra U$ is a graph $X(U_\ell)$ embedded in the Riemann surface $U(\C)$.  The vertices of $X(U_\ell)$ are the points of $U_\ell$ lying over $u$; we recall that these correspond to the $\ell^{2g}-1$ nontrivial $\ell$-torsion points on $A_u$.

On the other hand, we have a monodromy representation
\beq
\mu:  \pi_1(U(\C),u) \ra \Sp(H_1(A,\Z)) \cong \Sp_{2g}(\Z)
\eeq
The $\ell$-torsion in $A_u$ is canonically identified with the mod-$\ell$ homology group $H_1(A_u,\Z/\ell\Z)$.  So what happens if we start at a point $\tilde{u}$ in $U_\ell(\C)$ lying over $u$, and make $\tilde{u}$ ``follow its shadow" around a loop $\gamma_i$ in $U(\C)$?  The description of monodromy from the first section tells us that the path ends at $\mu(\gamma_i) \tilde{u}$.  In other words, the permutation on the preimages of $u$ induced by $\gamma_i$ is the same as that induced by the reduction mod $\ell$ of the symplectomorphism $\mu(\gamma_i)$.  This is precisely the same as saying that two vertices $\tilde{u}$ and $\tilde{u}'$ of $X(U_\ell)$ are adjacent precisely when $\tilde{u}'= \mu(\gamma_i) \tilde{u}$ for some $i$.  In other words, $X(U_\ell)$ is a quotient of the Cayley graph $X_\ell$ we described in section \ref{s:superstrong}.

This is the crucial point.  For a given family of abelian variety we may have no easy way of checking whether the monodromy group $\Gamma$ is thin or arithmetic.  But it doesn't matter!  The theorems of Helfgott, Pyber-Szab\"{o}, Breuillard-Green-Tao, and Golsefidy-Varju work equally well in either case to tell us that $\Gamma$ has superstrong approximation.  It is precisely this robustness that makes Theorem~\ref{th:sample} possible.

Superstrong approximation for $\Gamma$ tells us that the $X_\ell$ have a spectral gap; that is, the nontrivial eigenvalues of the {\em discrete} Laplacians on these graphs is bounded away from $0$ as $\ell$ grows.  (Hidden here is a certain amount of work necessary to pin down the possibilities for $G$, the Zariski closure of monodromy; this is necessary in order to ensure that superstrong approximation applies.  This is covered in \S 5 of \cite{ehk}.)

What is the relationship between the Laplacian of the graph and the geometry of the surface in which the graph is embedded?   It can be illustrative to visualize the case where all the $U_\ell$ have genus $0$; that is, they are homeomorphic to the sphere.  Now a graph embedded in a sphere is a planar graph, and it is known that a sequence of planar graphs with increasingly many vertices cannot have a spectral gap; intuitively, graphs with a spectral gap have very quickly dispersing random walks, while the random walk on a planar graph is constrained to linger for a while in the neighborhood of the plane where it starts out.

More generally, we should expect $X_\ell$ to serve as a kind of combinatorial approximation to the geometry of $U_\ell$.   Speaking roughly:  the graph is like a wire frame along which $U_\ell$ is stretched, and the graph separates $U_\ell$ into many small components, each one of which is a copy of a patch on $U$.  So the geometry of $U_\ell$ is, in a sense, a hybrid of the geometry of $X_\ell$ and the geometry of $U$.  As $\ell$ grows, and along with it the number of vertices of $X_\ell$, the graph plays a proportionately greater and greater role in the geometry of $U_\ell$ and the shape of $U$ a smaller and smaller one. 
  
This  intuition is beautifully confirmed and extended by results of Kelner, Brooks, and Burger, which together provide us with a ``comparison principle" between the Laplacian spectrum of the Riemann surface $U_\ell(\C)$ and the (discrete) Laplacian spectrum of the embedded graph $X_\ell$.  In particular, the combination of their theorems shows that the sequence of Riemann surfaces $U_\ell(\C)$ has a spectral gap whenever the sequence of graphs $X_\ell$ does.  This completes the proof.

\begin{rem}  For the arithmetic applications in \cite{ehk}, the full strength of superstrong approximation was not needed, and we used a weaker condition instead.  We denote sequence of graphs $\set{X_\ell}$ as {\em esperantist} when $\lambda_1(X_\ell)$ decays no more quickly than a power of the logarithm of $|X_\ell|$.  At the time \cite{ehk} was written, results guaranteeing superstrong approximation in the required generality were not available, but the theorems of Pyber and Szab\"{o}~\cite{pyberszabo} established the esperantist condition in the contexts we needed.  Though we now have superstrong approximation for thin groups in much greater generality, it still seems useful to keep in mind that in applications it is often sufficient to have at hand the easier condition on $\Gamma$ that the Cayley graphs attached to $\Gamma$ are esperantist (``slightly superstrong approximation...?")  On the other hand, the full strength of expansion allows one to conclude, not only that the gonality of $U_\ell$ grows without bound as $\ell \ra \infty$, but that the gonality is bounded below by a constant multiple of $\deg(U_\ell/U)$; in this case we say the gonality has {\em linear growth}.

\end{rem}

\begin{rem}  Many more applications of the basic method appear in \cite{ehk}.   For instance, if the monodromy group $\Gamma$ is Zariski dense in  $\Sp_{2g}(\Z)$, it follows from Theorem 4 of \cite{ehk} that there are only finitely many $t$ of degree at most $d$ over $k$ such that $A_t$ is isogenous to a product of abelian varieties of smaller dimension.
\end{rem}

\begin{rem}  The theorems in \cite{ehk} require that the families of abelian varieties be defined over fields of characteristic $0$, since everything rests on the ability to compare the geometry of Riemann surfaces with that of graphs drawn thereon.  The {\em statements} of the theorems, however, largely make sense in characteristic $p$; it seems very interesting to wonder whether the results continue to be true in this case, and in particular whether any analogue of the argument by superstrong approximation can be made to work.
\end{rem}

\begin{rem}  Work of Cadoret and Tamagawa, carried out independently around the same time as \cite{ehk}, gives lower bounds for the gonality of some of the covers discussed here, by very different methods (\cite{cado:2},\cite{cado:gonality}.)   In particular, their results are purely algebraic, and can be applied in characteristic $p$ as well as characteristic $0$.
\end{rem}

\section{Speculations:  non-abelian superstrong approximation}

We can phrase the approximation criteria of the previous sections in a different way.  Think of $\GL_n(\Z)$ as the automorphism group of the free abelian group $\Z^n$.  Then, for any finite abelian group $A$ which can be generated by $n$ elements, $\GL_n(\Z)$ acts by permutations on $\Epi(\Z^n,A)$, the set of surjections from $\Z^n$ to $A$:  we thus have a homomorphism
\beq
\Aut(\Z^n) \ra \Sym(\Epi(\Z^n,A))
\eeq

Strong approximation for a subgroup $\Gamma$ of $\GL_n(\Z)$  is essentially the statement that the restriction of this action to $\Gamma$  is {\em transitive} for all $A$; equivalently, the Cayley-Schreier graphs attached to this action are connected.  Superstrong approximation, on the other hand, implies the stronger statement that those graphs form an expander family.

It is natural to ask about non-abelian analogues of this story;  for instance, one can replace $\Z^n$ with a free group $F_n$ on $n$ generators, and replace the finite abelian group $A$ with a general finite group $G$.  It turns out that the analogues of strong and superstrong approximation have already appeared in various interesting contexts.

\subsection{The product replacement algorithm}

(Main reference:  Lubotzky-Pak~\cite{lubotzkypak}, Lubotzky~\cite{lubotzkysurvey})

If $G$ is a finite group, then $\Epi(F^n,G)$ is precisely the set of $n$-element generating sets $(g_1, \ldots, g_n)$ of $G$.  The {\em product replacement algorithm} of Leedham-Green and Soicher is a random walk on this set, which proceeds as follows:  choose distinct $i,j$ in $1..n$ at random,  and then multiply $g_i$ either on the left or right (chosen by coinflip) by $g_j$ or $g_j^{-1}$ (also chosen by coinflip.)

This random walk appears in practice to converge extremely quickly to the uniform distibution on generating $n$-tuples, making it an effective way to generate random elements in a group without having to enumerate the entire group.  In fact, the product replacement algorithm is known to converge to uniform distribution in polynomial time when $n$ is allowed to grow (slowly) with $|G|$, by a result of Pak~\cite{pak:prapoly}.  In the present paper we will concentrate mostly on the case where $n$ is fixed and $G$ changes.

The product replacement graph is in fact a Cayley-Schreier graph for the permutation action
\beq
A^+ \ra \Sym(\Epi(F^n,G))
\eeq
where $A^+$ is a subgroup of $\Aut(F_n)$ of index $2$, generated by automorphisms corresponding to the $4n(n-1)$ operations of the product replacement algorithm.  (We will casually fail to distinguish between $\Aut(F_2)$ and its finite-index subgroup $A^+$ in what follows.)  The rapid convergence of the random walk naturally suggests the question:  are these graphs expanders?  In other words, if $G_1,G_2,\ldots$ is a sequence of finite groups, we may ask:

\begin{question}
Does $\Aut(F^n)$ have superstrong approximation with respect to its image in $\Sym(\Epi(F^n,G_i))$?
\label{qu:superstrongpra}
\end{question}
The above question would certainly have a positive answer if $\Aut(F^n)$ had property $T$ (or even property $\tau$ for a suitable class of finite quotients.)   Under such a hypothesis, quantitative bounds for the performance of the product replacement algorithm are given in \cite{lubotzkypak}.  Whether $\Aut(F^n)$ has property $T$ in general is a well-known open question, about whose answer there is no clear consensus.  We know that $\Aut(F^n)$ does {\em not} have property $T$ when $n=2,3$~\cite{grunewaldlubotzky}.

One very interesting partial result was proved by Gamburd and Pak~\cite{gamburdpak}:  they show that Question~\ref{qu:superstrongpra} has a positive answer for the family of groups $\PSL_2(\F_p)$, under the hypothesis that $\PSL_2(\F_p)$ has {\em uniform expansion}; that is, there is a constant $\epsilon > 0$ such that, for all $p$, and all Cayley graphs of $\PSL_2(\F_p)$ with a specified number of generators, the spectral gap is at least $\epsilon$.  A very recent theorem of Breuillard and Gamburd shows that there is a family $P$ of primes of density $1$ such that $\PSL_2(\F_p)$ has uniform expansion as $p$ ranges over $P$, but the general case remains mysterious.

It's important to keep in mind that, for all our talk about superstrong approximation, even the prior question of strong approximation is poorly understood in the non-abelian case.  Recall that strong approximation has to do with the {\em transitivity} of the action on the set of epimorphisms.  There are many $(G,n)$ for which we know that $\Aut(F_n)$ acts transitively on $\Epi(F^n,G)$.  But there are many examples of $(G,n)$ such that the action is non-transitive; see \cite{pak:whatdoweknow} and \cite{lubotzkysurvey} for a full discussion of what is known and what is conjectured about this problem.

So far, no thin groups have appeared in this discussion.  It is not even precisely clear how the notion of ``thin subgroup" of $\Aut(F_n)$ should be defined.  Still, we can imagine what kind of questions might be asked in this setting.  For instance, one could ask about the following strengthening of Question~\ref{qu:superstrongpra}:

\begin{question}  Let $\Gamma$ be a finitely generated subgroup of $\Aut(F_n)$ whose image in $\Sym(\Epi(F^n,\PSL_2(\F_p))$ is the same as that of $\Aut(F_n)$ for all $p$.  Do the Cayley-Schreier graphs of the action of $\Gamma$ on $\Epi(F^n,\PSL_2(\F_p))$ form an expander family?
\end{question}

It goes without saying that it would be of interest to ask all of the above questions with $\Aut(F_n)$ replaced by the braid group on $n$ strands -- or, for that matter, the mapping class group of any oriented surface whose fundamental group is isomorphic to $F_n$.  This is more in the spirit of monodromy, since mapping class groups, unlike automorphism groups of free groups, are fundamental groups of natural moduli spaces -- namely, moduli spaces of punctured curves.

\subsection{Square-tiled surfaces and approximation for $\Aut(F_2)$}

The action of $\Aut(F_2)$ on $\Epi(F_2,G)$ also appears in the theory of {\em square-tiled surfaces}.  A square-tiled surface is just a connected degree-$d$ cover of a torus branched at a single point.  Such a cover can be thought of as a map from $\pi_{1,1}$ to $S_d$, where $\pi_{1,1}$ is the fundamental group of the once-punctured torus; it is a free group on $2$ generators.  

There is a natural obstruction to transitivity of the action of $\Aut(F_2)$ on $\Hom(F_2,S_d)$; namely, the conjugacy class of the puncture on the torus (alternately:  the commutator of a pair of free generators of $F_2$) is preserved up to conjugacy by $\Aut(F_2)$.  Since we have not demanded that our homomorphisms be surjective, we can break down $\Hom(F_2,S_d)$ further according to the image of the homomorphism.

Precisely:  Let $\mu$ be a conjugacy class of $S_d$, i.e. a partition of $d$, and let $G$ be a transitive permutation group of rank $d$ (i.e. a transitive subgroups of $S_d$ defined up to conjugacy.)  Then we denote by $Hom_\mu(F_2,S_d,G)$ the set of homomorphisms from $F_2$ to $S_d$ which send the puncture class to $\mu$ and whose image is $G$.  This set also carries an action of $S_d$, by conjugation on the right, and we denote the set of orbits of this action by $H^1_\mu(F_2,S_d,G)$.  

Now $\Aut(F_2)$ acts on $H^1_\mu(F_2,S_d,G)$, and in fact this action factors through its quotient $\Out(F_2)$.  We denote the associated Cayley-Schreier graph $X(d,\mu,G)$, and we can now ask questions of the same genre we have considered before:  for which $(d,\mu)$ is $X(d,\mu,G)$ connected?  For which sequences of $(d,\mu,G)$ do the $X(d,\mu,G)$ form an expander family?

We remark that $\Out(F_2)$ is naturally identified with $\GL_2(\Z)$ via its action on $F_2^{ab} = \Z^2$, and its index $2$ subgroup is in turn naturally identified with the mapping class group $\Gamma_{1,1}$ of the once-punctured torus, which is the fundamental group of $\MM_{1,1}$, the moduli space of genus-$1$ curves with a single marked point.  Under these identifications, the action of $\Gamma_{1,1}$ on the finite set $H^1_\mu(F_2,S_d,G)$ corresponds to a finite cover of $\MM_{1,1}$; this is precisely a moduli space $\HH_{d,\mu,G}$ of square-tiled surfaces, where the map $\HH_{d,\mu,G} \ra \MM_{1,1}$ sends a square-tiled surface to the punctured torus which it covers.

The moduli space $\HH_{d,\mu,G}$ is a {\em Hurwitz space} and is also an example of a {\em Teichm\"{u}ller curve}.  Their geometry and their relationship with moduli spaces of higher genus curves has been much studied -- see e.g. \cite{mcmullen:teich} and \cite{dawei}.  Any ``superstrong approximation" result for a family of $X(d,\mu,G)$ could in principle be used to give lower bounds for the gonality of $\HH_{d,\mu,G}$ along the lines explained in the previous section.  And of course one could go further and ask about approximation results for infinite-index subgroups of $\Gamma_{1,1}$ (though it is not clear at first glance whether such results have interesting implications in this context.)

We note that it is not possible that only finitely many $\HH_{d,\mu,G}$ have a component with low gonality, for the simple reason that every algebraic curve over a number field is birational to a connected component of $\HH_{d,\mu,G}$ for some $G$~\cite{eciatc}.  So any general statement of growing gonality would have to be restricted to a special family of $(d,\mu,G)$.

\section{Speculations:  superstrong approximation for Galois groups and the Bogomolov property}

In section~\ref{s:ehk} we explained how results on superstrong approximation could provide lower bounds for gonality in a family of covers of algebraic curves over characteristic $0$ fields.  The group to which the approximation results were applied was the image of the topological monodromy group of the base (a free group) in the arithmetic lattice $\Sp_{2g}(\Z)$.

If one tries to do the same thing over a more general base -- for instance, a base in characteristic $p$ -- one immediately encounters the problem that there is no longer a {\em discrete} fundamental group; all that remains is the \'{e}tale fundamental group, a topologically finitely generated profinite group.

What might superstrong approximation mean for a profinite, rather than a discrete group?  Let $\Pi$ be a profinite group which is topologically finitely generated.  If $\pi_1, \dots, \pi_k$ are topological generators, the subgroup $\Gamma$ they generate is dense in $\Pi$.  We say $\Gamma$ has property $\tau$ with respect to $\Pi$ if the Cayley graphs induced by $\Gamma$ on the finite quotients of $\Pi$ form an expander family.  We might wish to say that $\Pi$ satisfies superstrong approximation if the dense subgroup $\Gamma$ has property $\tau$ with respect to $\Pi$.  But there's a problem; it follows from a result of Alon, Lubotzky and Wigderson~\cite{alonlubotzkywigderson} show that there are profinite groups $\Pi$ and dense finitely generated subgroups $\Gamma,\Gamma' \subset \Pi$ such that $\Gamma$ has property $\tau$ with respect to $\Pi$, but $\Gamma'$ does not.  (Kassabov has shown~\cite{kassabov} that this can hold even when the profinite completions of $\Gamma$ and $\Gamma'$ are actually {\em isomorphic} to $\Pi$.)

One can instead place on $\Pi$ the stronger condition that {\em every} dense finitely generated subgroup $\Gamma$ has property $\tau$ with respect to $\Pi$.\footnote{There is a natural temptation to call this criterion ``superduperstrong approximation," but this terminology is perhaps better kept quarantined in a footnote.}  In this case we say that $\Pi$ has property $\hat{\tau}$.  This is very close to asking that the Cayley graphs attached to generating sets of size $k$ of finite quotients $Q$ of $\Pi$ satisfy a spectral gap which is independent of $Q$.  In other words, $\Pi$ needs to have {\em uniform expansion} as discussed in the previous section.  Uniform expansion is for the moment out of reach, but it is widely believed to hold, for example, when $\Pi = \Sp_{2g}(\Zhat)$, the relevant case for the current discussion.

Suppose $\Sp_{2g}(\Zhat)$ has property $\hat{\tau}$, let $U$ be a complex algebraic curve, and let $\rho:  \pi_1(U(\C),u) \ra \Sp_{2g}(\Zhat)$ be a representation of the fundamental group of $U$.  For each $\ell$, the reduction $\bar{\rho}_\ell: \pi_1(U(\C),u) \ra \Sp_{2g}(\Z/\ell\Z)$ gives rise to a cover $U_\ell$ of $U$, and the arguments of the previous section show that the $U_\ell$ grow in gonality.

One might ask the same question for a curve $U$ over a finite field $\F_q$.  We slightly modify the question in order to avoid introducing the machinery of the \'{e}tale fundamental group.  Suppose given, for an infinite sequence of primes $\ell$, an \'{e}tale cover $U_\ell \ra U$ such that $(U_\ell)_{\Fqbar} \ra U_{\Fqbar}$ is a Galois cover with Galois group $\Sp_{2g}(\Z/\ell\Z)$.  For instance, you might make such a cover by adjoining the $\ell$-torsion points of an abelian $g$-fold $A$ over $U$.

\begin{question}  Is the gonality of $U_\ell$ bounded below by a positive constant multiple of $[U_\ell:U]$?  What if we impose the extra condition that the $U_\ell$ are obtained by adjoining the $\ell$-torsion points on a fixed abelian $g$-fold $A/U$?
\label{qu:fqquestion}
\end{question}

We have focused on the case of curves over finite fields because of an intriguing connection with a different topic in number theory, which we now explain.  Write $K_\ell$ for the function field of $U_\ell$.  To say $U_\ell$ has gonality greater than $n$ is to say that there is no rational function $f$ on $U_\ell$ with at most $n$ poles and $n$ zeroes.   This is equivalent to the assertion that there is no non-constant element $f$ of $K_\ell^*$ whose (logarithmic) Weil height $h_{K_\ell}(f)$ is at most $n$.  The {\em absolute height} of an element $f$ of an extension of $K_\ell$ is just $h_{K_\ell}(f) / [U_\ell:U]$.  In other words, Question~\ref{qu:fqquestion} asks whether the {\em absolute} height of an element of $K_\ell^*$ is bounded below by a constant if it is nonzero.

This condition already has a name -- it is called the {\em Bogomolov property}:

\begin{defn}  Let $K$ be a global field.  A set $S$ of algebraic numbers over $K$ has the {\em Bogomolov property} if the set of positive absolute heights of elements of $S$ is bounded away from $0$.
\end{defn}

Note that this definition makes sense not only when $K$ is a function field (as in the case described above) but also when $K$ is a number field; in fact, it is the case $K=\Q$ that has been most energetically studied.  Some sets known to have the Bogomolov property over $\Q$ are the field of totally real algebraic numbers (Schinzel), the maximal abelian extension of $\Q$ (Amoroso-Zannier) and the maximal extension of $\Q$ totally split at a finite place $p$ (Bombieri-Zannier).  The recent paper of Amoroso, David, and Zannier~\cite{amorosodavidzannier} contains a good summary of what is known along with several striking new results.

Speculation piled atop speculation:  If we believe that $\Sp_{2g}(\F_\ell)$ has uniform expansion, and if we believe by analogy with the complex case that uniform expansion of Galois groups implies the Bogomolov property for the corresponding extensions of function fields, might this be the case for number fields as well?  There are many ways one could analogize Question~\ref{qu:fqquestion} to the number field case.  Here is a minimalist version:

\begin{question}  Let $A/\Q$ be an abelian variety and let $K_\ell$ be the number field $\Q(A[\ell])$.  Does the compositum of all the $K_\ell$ have the Bogomolov property?
\label{qu:minimal}
\end{question}

And here is a maximalist version:

\begin{question}  Let $F/\Q$ be an infinite Galois extension such that $\Gal(F\Q^{ab}/\Q^{ab})$ has property $\hat{\tau}$.  Does $F$ have the Bogomolov property?
\end{question}

A striking recent result of Habegger~\cite{habegger} shows that the answer to Question~\ref{qu:minimal} is affirmative when $A$ is an elliptic curve.  The theorems of \cite{ehk} described above can be used to show that the analogue of Question~\ref{qu:minimal} with $\Q$ replaced by $\C(t)$ has an affirmative answer as well; this requires lower bounds for the gonality, not only of the covers $U_\ell$, but of fiber products of multiple $U_{\ell_i}$ over $U$.  In this case, the relevant expansion theorem, which was not available when \cite{ehk} was written, is the theorem of Salehi Golsefidy and Varj\'{u}~\cite{golsefidyvarju}, which shows that expansion holds not only for the Cayley graphs associated to reduction modulo primes, but also modulo arbitrary squarefree integers.

As for the maximalist question, its interest is somewhat diminished by our lack of examples of groups known to have property $\hat{\tau}$.  This need not keep us from speculating.  For instance, $S$ be a finite set of primes, and let $K_S^t$ be the maximal extension of $\Q$ which is unramified outside $S$ and {\em tamely} ramified at $S$, and let $G_S^t$ be its Galois group.  Then $G_S^t$ has property {\em FAb}, which is to say that every finite-index subgroup of $G_S^t$ has finite abelianization.  This makes $G_S^t$ a reasonable candidate for having some form of spectral expansion property.  Caution is required, since $G_S^t$ is only conjectured to be topologically finitely generated.  One can circumvent this problem as follows.  For each integer $k$, we can consider the family of all Cayley graphs on finite quotients $G$ of $G_S^t$ endowed with subsets $(g_1, \ldots, g_k)$ generating $G$.  We say that $G_S^t$ has $\hat{\tau}_k$ if there is a uniform spectral gap for all such Cayley graphs.

\begin{question}  Does $G_S^t$ have property $\hat{\tau}_k$ for all $k$?  Does $K_S^t$ have the Bogomolov property?
\end{question}

Note that for $G_S^t$ to have property $\hat{\tau}_2$ would require that there are only finitely many $n$ such that $G_S^t$ has a quotient isomorphic to $S_n$.  It is not even known whether there are arbitrarily large $S_n$-extensions of $\Q$ unramified away from $S$; a recent preprint of Roberts and Venkatesh~\cite{robertsvenkatesh} argues that such extensions might indeed exist and proposes a construction, but at the same time suggests that these extensions are {\em not} likely to be tamely ramified at the primes in $S$.  Another hint that the first question is not completely unreasonable comes from Ershov's work in \cite{ershov} on Golod-Shafarevich groups with infinite quotients with property T.  But we should emphasize that we have no strong opinion as to the answer to either question.

\bibliographystyle{plain}
\bibliography{thingroups}

\end{document}